\documentclass[11pt]{amsart}
\usepackage{amssymb,amsthm,amscd}
\theoremstyle{remark}

\newtheorem{problem}{Problem}
\theoremstyle{plain}
\newtheorem{thm}{Theorem}
\newtheorem*{thm*}{Theorem}

\newtheorem{prop}{Proposition}
\newtheorem{cor}{Corollary}

\theoremstyle{definition}
\newtheorem{defin}{Definition}
\newtheorem{ex}{Example}
\def\paritem[#1]{\par\noindent\hbox to 2pc{\hss\rm#1}}
\def\sfrac#1#2{{\textstyle\frac{#1}{#2}}}
\newcommand\half{\sfrac 1 2}
\newcommand\ts{\mkern 2mu}
\DeclareMathOperator{\sspan}{\mathrm {span}}
\DeclareMathOperator{\Sec}{\mathrm {Sec}}
\DeclareMathOperator{\SU}{\mathrm {SU}}
\DeclareMathOperator{\U}{\mathrm U}
\DeclareMathOperator{\Ort}{\mathrm O}
\DeclareMathOperator{\SL}{\mathrm {SL}}
\DeclareMathOperator{\SO}{\mathrm {SO}}
\DeclareMathOperator{\CP}{\mathbb{C}\mathrm {P}}
\DeclareMathOperator\lrc{\text{\large\(\lrcorner\)}}
\DeclareMathOperator{\tbwedge}{\mbox{\large $\wedge $}}
\DeclareMathOperator{\tbotimes}{\mbox{\large $\otimes $}}
\DeclareMathOperator{\End}{\mathrm {End}}
\DeclareMathOperator{\id}{\mathrm {id}}
\DeclareMathOperator{\dir}{\mathrm {dir}}
\DeclareMathOperator{\rmd}{\mathrm d} %for d(
\DeclareMathOperator{\rd}{\mathrm d\!} %for dx
\DeclareMathOperator{\im}{\mathrm {Im}}
\DeclareMathOperator{\re}{\mathrm {Re}}
\newcommand{\rmi}{\mathrm i}
\newcommand{\sng}{{\em sng\/}}
\begin{document}
\title[robinson and hermite manifolds]{robinson manifolds
as the lorentzian analogs of hermite
manifolds}
\author{pawe\l\  nurowski and andrzej trautman}
\address{instytut fizyki teoretycznej, ho\.za 69, 00681 warszawa, poland}
\date{}
\thanks{Corresponding author: Andrzej Trautman.}
\thanks{The paper is in final form and will not be published elsewhere.}
%\thanks{Submitted for publication in {\it Differential Geometry and its
% Applications}.}
\begin{abstract}
A Lorentzian manifold is defined here as a smooth pseudo-Riemannian
 manifold with a metric tensor
of signature \((2n +1, 1)\). A Robinson manifold
is a Lorentzian manifold \(M\) of dimension \(\geqslant 4\)
with a subbundle \(N\)
of the complexification of \(TM\) such that the fibers
of \(N\to M\) are maximal totally null (isotropic) and \([\Sec N,
\Sec N]\subset \Sec N\). Robinson manifolds are close analogs
 of the proper Riemannian,
Hermite manifolds. In dimension 4, they correspond to space-times
of general relativity, foliated
by a family of null geodesics without shear. Such space-times,
 introduced in the 1950s
by Ivor Robinson, played an important role in the study of
solutions of Einstein's equations: plane and sphere-fronted waves,
the G\"odel universe,
the Kerr solution, and their generalizations, are among them. In this
survey article, the analogies between Hermite and Robinson manifolds are
presented in considerable detail.

\end{abstract}
\maketitle

\section{introduction and motivation from physics}
There is an interesting class of Lorentzian manifolds that  bear
a close analogy to the Hermite manifolds of proper Riemannian geometry.
They have been introduced and studied by physicists in the work on
 solutions of Einstein's equations, especially those
representing gravitational waves.  These {\it Robinson manifolds},
as we propose to call them, are little known to pure mathematicians.
This may be due, in part, to the fact that  physicists, in their
work,  used a local,
coordinate-dependent description of those manifold and did not
pay enough attention to
the geometrical motivation and interpretation of their results.
A good summary of this research by physicists is
in \cite{KSM}.

In this article, which is largely an expository survey, we describe the main
geometrical structures underlying Robinson manifolds and
emphasize their analogies with Hermite manifolds.

\subsection{Motivation from physics}
Let \(\mathbf E\) and \(\mathbf B\) be the vectors representing,
respectively, the
electric and magnetic fields in the Minkowski space-time \(\mathbb{R}^4\)
of special relativity theory.
Introducing \(\mathbf F=\mathbf E+\rmi\mathbf B\), one
can write Maxwell's equations in empty space in
the Riemann--Silberstein form (see
\cite{Silberstein1907} and p.~344 in \cite{Weber1910})
\begin{equation}\label{e:Max1}
    \rmi\frac{\partial}{\partial t}{\mathbf F}={\mathrm{curl}}\ts{\mathbf F}\quad\text{and}
\quad \mathrm{div}\ts{\mathbf F}=0.
\end{equation}
Among the solutions of \eqref{e:Max1} especially simple are the
{\it null\/} fields characterized by
\( \mathbf F^2=0\).
The property of \(\mathbf F\) to be null can be linearized: it is equivalent to
the statement
\begin{equation}\label{e:Fnull1}
  \text{there exists a unit vector \(\mathbf n\) such that \(
\mathbf n\times\mathbf F=\rmi \mathbf F\).}
\end{equation}
Introducing an orientation in \(\mathbb R^4\) defined by the
 form \(\rd t\wedge\rd x\wedge\rd y\wedge\rd z\)
so that Hodge duality of 2-forms is given by
\begin{equation*}
\star(\rd t\wedge\rd x)=\rd y\wedge\rd z,\quad
\star(\rd y\wedge\rd z)=-\rd t\wedge\rd x,\quad\text{etc.,}
\end{equation*}
putting
\begin{equation*}F=F_x(\rd t\wedge\rd x-\rmi \rd y\wedge\rd z)+
\mathrm{cycl.}\quad
\text{and}\quad \kappa=\rd t-n_x\rd x-n_y\rd y-n_z\rd z,
\end{equation*}
 one has
\begin{equation}\label{e:starF}
    \star F=\rmi F
\end{equation}
and can write
\eqref{e:Max1} and \eqref{e:Fnull1} in the equivalent form
\begin{equation}\label{e:Max2}
    \rd  F=0,
\end{equation}
and
\begin{equation}\label{e:Fnull2}
 \text{there exists a 1-form \(\kappa\neq 0\) such that
\(\kappa\wedge  F=0\),}
\end{equation}
respectively.

The virtue of conditions \eqref{e:starF}-\eqref{e:Fnull2} is that, without
change of form, they are meaningful on every oriented,
4-dimensional  Lorentzian manifold \((M,g)\). (In fact, conformal geometry
of Lorentzian signature is enough and one can generalize to
a \(2n\)-dimensional manifold by assuming, in addition,
 that \(F\) is a decomposable
\(n\)-form.) A  4-dimensional  Robinson
manifold can be provisionally defined as a Lorentzian manifold admitting
a nowhere zero, complex-valued 2-form \(F\) such that conditions
\eqref{e:starF}-\eqref{e:Fnull2} hold.
The vector field \(k\)
associated by \(g\) with \(\kappa\) is
{\it null\/}. (Pure mathematicians say:
{\it isotropic}, but this is a misnomer.
The term isotropic was introduced, in this context, by Ribaucour
(see Ch. 4 in
 \cite{Klein})
in the study of  complex
Euclidean geometry: if
\(\mathbb C^2\) is endowed with the quadratic form \((z_1,z_2)
\mapsto z_1^2+z_2^2\), then a rotation by the angle \(\alpha\)
transforms the vector \((1,\rmi )\) into
\((\exp \rmi\alpha ,\rmi \exp \rmi\alpha )\). This vector is
isotropic in the sense that its direction does not change under
rotations. But null directions in higher dimensions are not
invariant under rotations. Cartan   had the good idea of calling such
directions in \(\mathbb R^4\) {\it
optical\/}, but this name has not caught on.)  The field \(k\)
 defines a foliation (physicists say:
congruence) of \(M\) by null geodesics (Mariot's theorem; see
\cite{RT86} and the references given there).
 Ivor Robinson \cite{Rob61} found a
necessary condition on the foliation, which
is also sufficient in the analytic case,
but not otherwise \cite{Taf85}, for the existence of a nowhere vanishing
solution \(F\)  of \eqref{e:starF}-\eqref{e:Fnull2}. In the
physicists' language this condition is expressed by saying that \(k\)
should generate a {\it shear-free  null geodetic\/} (\sng)
congruence;
see \S \ref{ss:Rob}.

\subsection{Historical remarks  and plan of the article}
In 1910, Harry Bateman \cite{Bat10} discovered a class of
transformations, more general than conformal changes of the metric,
that can be used to transform null solutions of Maxwell's equations
into similar solutions; this work can be considered to be a
precursor of the `optical' ideas we are describing here; see
\cite{RT3,Trautman99} and Theorem \ref{t:Bat}.
In a short note of 1922, \'Elie Cartan \cite{C1} mentioned the existence of four
principal  optical (null) directions associated with a
non-conformally flat Lorentz  4-manifold.
 He also pointed out that, in the case of
the Schwarzschild space-time, these directions degenerate to form two
pairs of double optical directions. Cartan's observations went
unnoticed for almost 50 years. In the meantime and independently,
A.Z.~Petrov \cite{Petrov} devised an algebraic classification of the
Weyl tensor (of conformal curvature)
of a Lorentzian manifold and F.A.E.~Pirani \cite{Pirani} clarified
its physical significance. Using Weyl (two-component) spinors,
Roger Penrose \cite{Pen60} sharpened the Petrov classification
and gave a new derivation of the four null directions;
this is recalled here
in \S \ref{sss:aclass}. This and
subsequent work by Penrose (see \cite{PenRin} and the references
given there) has had a decisive influence on the development of
the subject. From the perspective of this article,  most significant
was the discovery by I.~Robinson \cite{Rob61} of the shear-free
property of congruences of null geodesics and their relation to
null electromagnetic fields (\S \ref{ss:Rob}).
To make the article self-contained and moderately complete, we
have included several classical theorems related to its subject,
with references to literature instead of proofs.
 In particular, in Section \ref{s:GS}
we present the Goldberg--Sachs theorem on the connection
between the existence of \sng\
 congruences  and the degeneracy of the principal null
directions in Einstein manifolds, as well as its generalization to
the proper Riemannian case.
A theorem due to R.P.~Kerr, giving all \sng\ congruences
 in  Minkowski space-time is presented
 in considerable detail in \S
\ref{sss:Kerr} and \ref{ss:Kerrrev}. In the last Section,
 we briefly describe twistor
bundles, an important concept that emerged in connection with the study
of \sng\ congruences. There is a wealth of literature
 on Penrose's twistor ideas, in both the Lorentz and proper
 Riemannian cases \cite{AtH78,Pen67,PenMcC73,PenRin,Wells79}.
 Recent surveys are in \cite{HMTTW98}.

\section{notation and terminology}
Our notation and terminology are essentially standard; see, e.g.,
\cite{Besse,KobNom,LawsonMich89}.
The exterior algebra associated with a vector space
 \(W\) is \(\tbwedge W\); the symbols \(\otimes\), \(\wedge\) and \(\lrc\)
 denote the tensor, exterior and  interior products, respectively.
 We use the Einstein summation convention over repeated indices.
The canonical map of \(W\smallsetminus\{0\}\) onto the associated
projective space \(\mathrm{P}(W)\) is denoted by \(\dir\) and
we write \(\CP_n\) for \(\mathrm P(\mathbb{C}^{n+1})\).
 A {\it quadratic space\/} is defined as a pair \( (V,g)\), where \(
 V\) is a  finite-dimensional vector space
 over \(\Bbbk=\mathbb{R}\) or \(\mathbb{C}\),
  and \( g:V\to V^*
 \) is a symmetric (\(g^*=g\)) isomorphism. To save on
 notation, we use the same letter \(g\) for the {\it metric tensor\/}
  \(g\in V^*\otimes_{\mathrm{sym}}V^*\)  associated with that isomorphism
 so that \(g(u,v)=\langle u, g(v)\rangle\) and
  \(v\mapsto g(v,v)\) is a quadratic form.
  For the symmetrized tensor product of 1-forms we use
  the notation of classical differential geometry,
  i.e., if \(\alpha,\beta\in V^*\), then
  \(2\alpha\beta=\alpha\otimes\beta+\beta\otimes\alpha\).
  This convention allows us to write the metric tensor as
   \(g=g(e_\nu)e^\nu=g_{\mu\nu}e^\mu e^\nu\), where \((e^\mu)\) is the coframe
   dual to \((e_\mu)\) and \(g_{\mu\nu}=g(e_\mu,e_\nu)\).
If \(N\subset V\), then \(N^\perp\) is the set of all elements of
\(V\) orthogonal to every  element of \(N\).
  The Hodge dual of \(\alpha\) is denoted by
\(\star\alpha\).

 All manifolds and maps among them are assumed to be smooth
(of class \(C^\infty\)) or real-analytic.
 Manifolds are finite-dimensional, but not necessarily compact.  If \(
 f:M'\to M \) is a map of manifolds, then \( Tf:TM'\to TM \) is the
 corresponding tangent (derived) map and
 \( T_{x}M\subset TM\) is the tangent vector space to \( M\) at \(
 x\). The map \(f\) is an
 immersion (resp., submersion) if \(Tf\), restricted to every tangent
 vector space, is injective (resp., surjective);
 an injective immersion is an
 embedding and defines \(M'\) as a submanifold of \(M\).
   If \( \pi:E\to
 M\) is a fiber bundle over a manifold \( M\), then \( E_{p}=
 \pi^{-1}(p)\subset E\) is the fiber over \( p\in M\).
 A map \(f:M'\to M \) gives rise to the induced bundle
  \(f^{-1}E\to M'\) such that \((f^{-1}E)_p=E_{f(p)}\) for every
  \(p\in M'\). If \(f\) is
an  immersion, then
 \(TM'\) is a subbundle of  \(f^{-1}TM\). The zero bundle is denoted by
 \(\pmb 0\).
A {\it Riemannian manifold\/}
 \( M \) is assumed to be connected; it has
 a metric tensor field \( g
 \) which is non-degenerate, but not necessarily definite; if it is,
 then \( (M,g) \) is said to be {\it proper\/} Riemannian.
 A {\it space-time\/} is a 4-dimensional manifold with a metric
 tensor of signature \((3,1)\).

 The module over \(C^\infty(M)\) of all sections
 of the vector bundle \(E\to M\) is denoted by \(\Sec E\).
If \(X\in\Sec TM\), then \(L(X)\) is the Lie derivative with
respect to \(X\). If \(\alpha\) is a differential form
on \(M\) and \(f:M'\to M\),
then \(L(X)\alpha=X\lrc\rd\alpha+\rmd(X\lrc\alpha)\) and \(f^*\alpha\)
is the pull-back of \(\alpha\) to \(M'\).
We abbreviate \(\partial/\partial x\)
to \(\partial_x\).
In Section \ref{s:CR} we summarize
 the definitions and notions related to CR structures needed in
 this paper; further details can be found in \cite{Jac90}.

 To save on  notation, we sometimes use the same letter to denote
 a vector space \(N\) with some structure and a fiber
 bundle \(N\to M\) with
 fibers carrying the same structure. Local sections of \(N\to M\)
 may be denoted by the same letters as elements of the vector space \(N\).

\section{algebraic preliminaries}
\subsection{Maximal, totally null subspaces of vector
spaces}
\noindent Consider a  complex quadratic space \((V,g)\).
Recall that a vector subspace \(N\) of \(V\) is said to be
{\it null\/} if \(N^\perp\cap N\neq \varnothing\) and
{\it totally null\/}
if \(N\subset N^\perp\). Assume now \(\dim V={2n}\);
if \(N\subset V\) is {\it maximal totally null\/} ({\it mtn\/}), then
 \(N^\perp=N\) so that \(\dim N=n\).
An orientation having been fixed, the Hodge duality map
 \(\star:\tbwedge V\to\tbwedge V \) can be
defined so that \(\star^2=\id\). If \((m_1,\dots,m_n)\) is
a frame in an {\it mtn\/} subspace \(N\), then
\begin{equation}\label{e:star}
\star(m_1\wedge\dots\wedge m_n)=\pm m_1\wedge\dots\wedge m_n.
\end{equation}
The {\it annihilator\/} of \(N\),
\begin{equation*}
  N^0=\{\mu\in V^*\mid \langle m,\mu\rangle=0\;\;\text{for every}\;\;
m\in N\}
\end{equation*}
is an {\it mtn\/} subspace of \(V^*\).
 The set of all {\it mtn\/} subspaces of a complex, \(2n\)-dimensional
 vector space has the structure of a complex manifold,
 diffeomorphic to the symmetric space \(\Ort_{2n}/\U_n\);
 its two connected components correspond to the two signs in \eqref{e:star}
  characterizing the {\it mtn\/} subspaces of positive and negative
  chiralities, respectively.

Let now \((V,g)\) be a Euclidean quadratic space, i.e. a real quadratic
space
such that the  form associated with \(g\) is positive-definite.
Assume that \(V\) is of positive even dimension.
An {\it mtn\/} subspace
\(N\) of the complexification \(W=\mathbb{C}\otimes V\)
defines a complex orthogonal structure \(J\) on \((V,g)\):
 this is so because
\(N\cap \bar{N}=\{0\}\) and one can put
\begin{equation}\label{e:defJ}
J(v)=\rmi v\quad\text{and}\quad J(\bar{v})=-\rmi \bar{v}
\quad\text{for}\quad v\in N.
\end{equation}
 Conversely, an orthogonal complex structure
\(J\) on \((V,g)\) defines the {\it mtn\/} subspace
\(N=\{v\in W\mid J(v)=\rmi v\}\).

Consider now    a {\it Lorentz   space\/} \((V,g)\), defined as
a real quadratic space
such that the quadratic form associated with \(g\) is
 of signature \((2n+1,1)\), \(n=1,2,\dots\).
 Let \(N\subset W=\mathbb{C}\otimes V\)
 be an {\it mtn\/} subspace. The intersection
 \(N\cap\bar{N}\) is the complexification of a null real line
 \(K\subset V\) and
\(
  N+\bar{N}=\mathbb{C}\otimes K^\perp
\).
There is a real null line \(L\) such that \(V=K^\perp\oplus L\).
The quotient \(K^\perp/K\) inherits from \((V,g)\) the structure
of a Euclidean quadratic space of dimension \(2n\) and there
is an orthogonal complex structure \(J\) on \(K^\perp/K\),
defined by
\(
  J(v \bmod \mathbb{C}\otimes K)=\rmi v \bmod \mathbb{C}\otimes K\)
for every \(v\in \mathbb{C}\otimes K^\perp
\).
Similarly, \(N^0\cap \bar{N}^0\) is the complexification of a real
null line and there is the isomorphism
\begin{equation}\label{e:kkappa}
  g:K\to \re N^0\cap \bar{N}^0
\end{equation}
obtained by restricting \(g:V\to V^*\) to \(K\).

\subsection{Spinor algebra in dimension 4}
Spinor calculus in dimension 4  provides
an economical, convenient description of many aspects of
the geometry of Riemannian manifolds of this dimension
\cite{LawsonMich89,PenRin}. Since there are so many
exhaustive presentations of this subject, it suffices
to give here the rudiments of spinor algebra in a form
adapted to our purposes.

If the  dimension of the real vector space
\(V\) is 4, then the complex vector space \(S\) of Dirac spinors is
also four-dimensional. Let \((e_\mu)\) be an orthonormal
frame in \(V\). A representation \(\gamma\)
of the Clifford algebra associated with \((V,g)\) in \(S\)
is given by the `Dirac matrices' \(\gamma_\mu=\gamma(e_\mu)\).
The endomorphism \(\gamma_5=\gamma_1\gamma_2\gamma_3\gamma_4\)
anticommutes with the Dirac matrices and \(\gamma_5^2=\id\)
if \((V,g)\) is  Euclidean and \(\gamma_5^2=-\id\) if
\((V,g)\) is Lorentzian. Putting \(\varGamma=\gamma_5\) in the first
and \(\varGamma=\rmi\gamma_5\) in the second case, one has
\(\varGamma^2=\id\).

 The spaces of `chiral' or Weyl spinors
are defined by
\begin{equation*}
  S_\pm=\{\varphi\in S\mid \varGamma\varphi=\pm\varphi\}.
\end{equation*}
Let \(W=\mathbb{C}\otimes V\) and, for \(v_1,v_2\in V\),
 put \(\gamma(v_1+\rmi v_2)=
\gamma(v_1)+\rmi\gamma(v_2)\), then \(\gamma(w)^2=g(w,w)\id\)
for every \(w\in W\). If \(\varphi\in S_\pm\) and \(\varphi\neq 0\),
then
\begin{equation}\label{e:Nvar}
  N(\varphi)=\{w\in W\mid \gamma(w)\varphi=0\}
\end{equation}
is an {\it mtn\/} subspace of \(W\) of the same chirality as \(\varphi\).

The transposed endomorphisms \(\gamma^*_\mu\) define
the contragredient representation of the Clifford algebra in
\(S^*\), which is equivalent to \(\gamma\): there is the
isomorphism \(B:S\to S^*\) such that \(\gamma^*_\mu=B\gamma_\mu B^{-1}\)
for \(\mu=1,\dots,5\).
 \(B\) restricts
to a symplectic form \(\varepsilon\)
 on each of the  spaces of
Weyl spinors \(S_+\) and \(S_-\). If \((e_A)\), \(A=1,2\), is a frame
in \(S_+\) and \((e^A)\) is the dual frame in \(S^*_+\), then
\(\varepsilon(e_A)=\varepsilon_{AB}e^B\).
The complex conjugate representation given by \(\bar{\gamma}_\mu\)
is also equivalent to \(\gamma\): there is an isomorphism
\(C:S\to \bar{S}\) such that \(\bar{\gamma}_\mu=C\gamma_\mu C^{-1}\)
and \(C\bar{C}=-\id\) in the Euclidean case and \(C\bar{C}=\id\)
for signature \((3,1)\).
The spinor \(\varphi_c=C^{-1}\bar{\varphi}\) is said (by physicists) to be the
{\it charge conjugate\/} of \(\varphi\in S\).

\subsection{The algebraic classification of Weyl
tensors}\label{sss:aclass}

The spaces \(S^4_+=\tbotimes^{4}_{\mathrm{sym}}S^*_+\)
and \(S^4_-=\tbotimes^{4}_{\mathrm{sym}}S^*_-\) are isomorphic
to spaces of tensors of rank \(4\)
over \(W=\mathbb{C}^4\), with symmetries of self-dual and
anti-self-dual Weyl (conformal curvature) tensors, denoted by
\(\mathsf C_+\) and \(\mathsf C_-\), respectively.
Consider  \(0\neq\psi\in S^4_+\):
 there is a frame \((e_A)\), \(A=1,2\), in \(S_+\)
such that the component
 \(\psi_{1\dots 1}=\psi(e_1,\dots,e_1)\)  is not zero.
 Given such a frame, let
\(\varphi(z)=ze_1+e_2\in S_+\), \(z\in\mathbb{C}\),
and consider the
complex polynomial \(p_\psi\) of degree \(4\),
\begin{equation*}
p_\psi(z)=\psi(\varphi(z),
\dots, \varphi(z))=\psi_{1\dots 1}z^4+\dots+\psi_{2\dots 2}.
\end{equation*}
 Let \(\{z_1,\dots,z_4\}\) be the set
of all roots of this polynomial; a root of multiplicity \(s\)
appears \(s\) times in the set.
Then
\begin{equation*}
  \psi=\psi_{1\dots 1}\varphi^1\otimes_{\mathrm{sym}}
   \dots\otimes_{\mathrm{sym}}\varphi^4,
  \quad\text{where}\quad \varphi_A^i=\varepsilon_{AB}\varphi(z_i)^B,\;\;i=1,\dots,4.
\end{equation*}
The spinors \(\varphi^i\) are {\it eigenspinors\/} (with eigenvalue 0)
of \(\psi\).
The {\it algebraic type\/} of \(\psi\) is the sequence
\({\pmb [}s_1\dots s_k{\pmb ]}\), \(1\leqslant s_1\leqslant \dots\leqslant s_k
\leqslant 4\), \(s_1+\dots+s_k=4\),
 of the multiplicities of the roots of \(p_\psi\).
In the generic case, all roots are simple, \(s_1=\dots=s_4=1\).
Otherwise, one says that \(\psi\) is {\it algebraically
degenerate\/}. An eigenspinor is said to be {\it repeated\/}
if its multiplicity \(s\) is larger than 1.

The enumeration of the possible degeneracies can be traced
back to Cartan \cite{C1}; physicists use it now  in a form due
to Penrose \cite{Pen60}:
\smallskip

\parbox{.6\linewidth}{

\paritem[(i)] Type I (non-degenerate) \({\pmb [}\ts 1111\ts{\pmb ]}\),

\paritem[(ii)] Type II \({\pmb [}\ts 112\ts{\pmb ]}\),

\paritem[(iii)] Type III \({\pmb [}\ts 13\ts{\pmb ]}\),

\paritem[(iv)] Type D (`degenerate') \({\pmb [}\ts 22\ts{\pmb ]}\),

\paritem[(v)] Type N (`null') \( {\pmb [}\ts 4\ts{\pmb ]}\).}

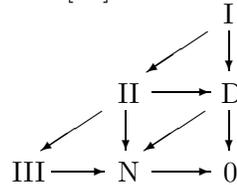
\parbox{.4\linewidth}{\begin{picture}(10,10)(-160,15)

\put(150,85){I}\put(152,82){\vector(0,-1){15}}
\put(110,55){II}\put(113,52){\vector(0,-1){15}}
\put(149,55){D}\put(152,52){\vector(0,-1){15}}
\put(70,25){III}\put(110,25){N}\put(150,25){0}
\put(85,29){\vector(1,0){20}}
\put(123,29){\vector(1,0){22}}
\put(123,59){\vector(1,0){22}}
\put(144,82){\vector(-3,-2){23}}
\put(104,52){\vector(-3,-2){23}}
\put(143,52){\vector(-3,-2){23}}
\end{picture}
}

%\smallskip
\noindent The \(0\) in the Penrose diagram
above represents a vanishing \(\psi\).
The arrows point towards more special cases. This classification of
complex, self-dual Weyl tensors is often associated with the name
of Petrov, who, however, recognized only three types (I, II and III).
The Weyl tensor of a complex Riemannian manifold decomposes into
its self-dual and anti-self-dual parts; their algebraic types are
independent.

In the case of real manifolds, one has to consider
separately each signature. We restrict ourselves to the proper
Riemannian and Lorentzian cases.

1. In the proper Riemannian case, the Weyl tensor decomposes
into the real, self-dual and anti-self-dual parts; they are
independent. The self-dual
part is represented by a spinor \(\psi\in S^4_+\) that satisfies a
suitable reality
condition which
implies  that the eigenspinors of
 \(\psi\) occur in pairs \((\varphi,\varphi_c)\).
Therefore, there are only two types of \(\psi\neq 0\): either
these two pairs are distinct (type I) or they coincide (type D).
Similar remarks apply to the anti-self-dual part of the Weyl tensor.
Therefore, the complete algebraic classification of the Weyl tensor
of a proper Riemannian 4-dimensional manifold contains 9 cases;
(I,I) is the most general case and (0,0) represents
conformally flat manifolds. The cases \((*,0)\) and \((0,*)\)
are referred to as  self-dual and anti-self-dual, respectively.

2. In the Lorentzian case, the real Weyl tensor
decomposes into its self- and anti-self-dual parts, which are complex,
 \(\mathsf C={\mathsf C}_+ +\mathsf C_-\), where
 \(\star{\mathsf C}_\pm=\pm\rmi{\mathsf C}_\pm\) so
 that  \(\bar{\mathsf C}_+=\mathsf C_-\). Therefore, the classification
is given by that of the complex, self-dual Weyl tensor presented above.
\section{cauchy--riemann manifolds}\label{s:CR}
\subsection{Almost CR manifolds}
\begin{defin}
An {\it almost Cauchy--Riemann manifold}  \(\mathcal M\)
of dimension \(2n+1\) is defined as a manifold with a distinguished subbundle
\(\mathcal{N}\) of \(\mathbb{C}\otimes T\mathcal{M}\), with
fibers of complex dimension \(n\),  such that
\(\bar{\mathcal{N}}\cap\mathcal{N}=\pmb 0\).
\end{defin}
One also says that \(\mathcal{M}\) has an almost CR structure.
The direct sum \(\bar{\mathcal{N}}\oplus\mathcal{N}\) is the complexification
of a bundle \(\mathcal{H}\subset T\mathcal{M}\) with \(2n\)-dimensional
fibers, endowed with \(J\in\Sec\End \mathcal H\)
 such that \(J^2=-\id_{\mathcal H}\); namely,
\(J(w+\bar{w})=\rmi (w-\bar{w})\) for every \(w\in\mathcal{N}\).

The annihilator \(\mathcal{N}^0\subset\mathbb{C}\otimes T^*\mathcal{M}\)
has fibers of complex dimension \(n+1\) and \(\overline{\mathcal{N}^0}
\cap\mathcal{N}^0\) is the complexification of a real line bundle.
The {\it canonical bundle\/} \cite{Jac87} of the almost CR structure,
\(
  \varOmega=\tbwedge^{n+1}\mathcal{N}^0
\),
is a complex line bundle over \(\mathcal{M}\) and
\begin{equation*}
  \mathcal{N}_p=\{w\in \mathbb{C}\otimes T_p\mathcal{M}\mid
  w\lrc \omega=0,\,\,0\neq\omega\in\varOmega_p,\,\,p\in\mathcal{M}\}.
\end{equation*}
There is a convenient, equivalent description of an almost CR structure
by an atlas of CR compatible charts: every point of \(\mathcal{M}\) has
a neighborhood \(\mathcal{U}\) admitting
 a collection of
1-forms
\begin{equation}\label{e:kappamu}
(\kappa,\mu^1,\dots,\mu^n)\;\text{with \(\kappa\) real and
\(\kappa\wedge\mu^1\wedge\dots\wedge\mu^n\wedge\overline{\mu^1}\wedge
\dots\wedge\overline{\mu^n}\neq 0\)}
\end{equation}
such that
\begin{equation}\label{e:NN0}
  \mathcal{N}^0_p=\sspan_p\{\kappa,\mu^1,\dots,\mu^n\}\quad\text{for every}
\quad p\in\mathcal{U}.
\end{equation}
The pair
\begin{equation}\label{e:CRchart}
(\mathcal{U},(\kappa,\mu^1,\dots,\mu^n))
\end{equation}
is a {\it CR chart\/}.
Given any other CR chart \((\mathcal{U'},(\kappa',{\mu'}^1,\dots,{\mu'}^n))\), on
the overlap \(\mathcal{U}\cap\mathcal{U'}\) one has
\begin{equation}\label{e:eqCR}
  \kappa'=a\kappa,\quad
 {\mu'}^\alpha=b^{\ts\alpha}\kappa+{b^{\ts\alpha}}_\beta\mu^\beta,
 \quad \alpha,\beta=1,\dots,n,
\end{equation}
where \(a\) is a real function, the \(b\ts\)s are complex and \(a\det b\neq 0\),
where \(b=({b^{\ts\alpha}}_\beta)\).
An almost CR manifold can be defined as an odd-dimensional
manifold with an atlas of compatible CR charts, their
compatibility being defined by \eqref{e:eqCR}.
The \((n+1)\)-form
\begin{equation}\label{e:defomega}
  \omega=\kappa\wedge\mu^1\wedge\dots\wedge\mu^n,
\end{equation}
is a  nowhere vanishing local section of \(\varOmega\to\mathcal{M}\)
defined on \(\mathcal{U}\).

Given \eqref{e:kappamu}, one puts
\begin{equation*}
  \rd \kappa=\rmi h_{\alpha\beta}\mu^\alpha\wedge\bar{\mu}^\beta+\dots,
  \end{equation*}
  where the dots stand for exterior products of pairs of the local basis
  1-forms other than the products
  \(\mu^\alpha\wedge\bar{\mu}^\beta\), \(1\leqslant \alpha,\beta
  \leqslant n\). The transformation \eqref{e:eqCR} induces the change
\begin{equation*}
  h'_{\alpha\beta}=ah_{\gamma\delta}{c^{\ts\gamma}}_\alpha
  {\bar{c}^{\ts\delta}}_\beta, \quad 1\leqslant
 \alpha,\beta,\gamma,\delta\leqslant n,
\end{equation*}
where \(c=({c^{\ts\alpha}}_\beta)\) is the inverse of the matrix \(b\).
The matrix \(h=(h_{\alpha\beta})\) is Hermitean and the signature
of the associated Hermitean {\it Levi form\/} is well-defined:
it does not change under the replacement \eqref{e:eqCR}. The almost
CR structure is said to be {\it non-degenerate\/} if \(\det h\neq 0\);
it is called {\it pseudo-convex\/} (sometimes: strongly pseudo-convex)
 if the associated Hermitean form is definite.

If the distribution \(\ker\kappa=\mathcal{H}\) is integrable,
\(\kappa\wedge\rd\kappa=0\), then the CR structure is said to
be {\it trivial\/} and, locally, \(\mathcal M=\mathbb{R}\times\mathbb{C}^n\).
In dimension three, non-triviality of a CR structure is equivalent to
its pseudo-convexity.
\subsection{CR manifolds}
\begin{defin}
A {\it Cauchy--Riemann manifold\/} \((\mathcal{M},\mathcal{N})\)
is an almost CR manifold
characterized by the bundle \(\mathcal{N}\to \mathcal{M}\), satisfying
the integrability condition
\(
  [\Sec \mathcal{N},\Sec\mathcal{N}]\subset\Sec\mathcal{N}
\).
\end{defin}
The integrability condition is equivalent to
\begin{equation*}
  \rd\Sec \mathcal{N}^0\subset\Sec\mathcal{N}^0\wedge\Sec(\mathbb{C}
  \otimes T^*\mathcal{M}).
\end{equation*}
In terms of a CR chart  \eqref{e:CRchart}
of \(\Sec \mathcal{N}^0\) this is
  equivalent to
\begin{equation}\label{e:int}
  \rd\kappa\wedge\omega=0\quad\text{and}
  \quad \rd\mu^\alpha\wedge\omega=0\quad\text{for}
  \quad \alpha=1,\dots,n.
\end{equation}
Clearly, every 3-dimensional almost CR manifold
is a CR manifold; we refer to it as a {\it CR space}.

If the canonical bundle \(\varOmega\) admits, for every
\(\mathcal U\) in the atlas, a closed local section \(\omega\) nowhere zero on
\(\mathcal U\),
then the integrability conditions \eqref{e:int} follow
from \(\kappa\wedge\omega=0\) and \(\mu^\alpha\wedge\omega=0\),
 \(\alpha=1,\dots,n\).

The chart \eqref{e:CRchart} is said to be locally
{\it embeddable\/} (sometimes: realizable)
if the {\it tangential CR equation}
\begin{equation}\label{e:embed}
  \rd z\wedge\omega=0
\end{equation}
has \(n+1\) solutions \(z_1,\dots,z_{n+1}\) such that
\begin{equation*}
  \sspan_p\{\rd z_1,\dots,\rd z_{n+1},\rd \bar{z}_1,\dots,
  \rd\bar{z}_{n+1}\}=\mathbb{C}\otimes T_p^*\mathcal{M}
  \quad\text{for every \(p\in\mathcal{U}\)}.
\end{equation*}
One then has the exact local section \(\omega=\rd z_1\wedge\dots\wedge
\rd z_{n+1}\) of the canonical bundle and
the map \(z:\mathcal{U}\to\mathbb{C}^{n+1}\approx\mathbb{R}^{2n+2}\),
 \(z=(z_1,\dots,z_{n+1})\),
is  an immersion. A CR manifold is locally embeddable if it has
a CR atlas of locally embeddable charts. Every analytic CR manifold
is locally embeddable \cite{AndrHill72}.

Let \(\mathcal M\) be now an  embeddable CR space
so that  there are two solutions \(z_1\) and \(z_2\) of \eqref{e:embed}
and a real-valued smooth function \(G\) on \(\mathbb{C}^2\) such
that
\begin{equation}\label{e:G=0}
G(z_1,z_2,\bar{z}_1,\bar{z}_2)=0\quad \text{and}\quad \rd G\neq 0.
\end{equation}
  One can then take
\begin{equation}\label{e:emb}
   \kappa=\rmi  (\frac{\partial G}{\partial z_1}\rd z_1+
   \frac{\partial G}{\partial z_2}\rd z_2),\quad \mu=
\overline{\frac{\partial G}{\partial z_2}}\rd z_1
-\overline{\frac{\partial G}{\partial z_1}}\rd z_2.
\end{equation}

\subsection{CR submanifolds}
\begin{defin}\label{d:sub}
Let \((\mathcal{M},\mathcal{N})\) and
\((\mathcal{M'},\mathcal{N}')\)  be  CR manifolds
 of dimension \(2n+1\) and \(2n-1\), respectively.
If  \(\mathcal{M'}\) is a submanifold of \(\mathcal{M}\)
 with an embedding \(f:\mathcal{M'}\to\mathcal{M}\) and
 \(\mathcal{N'}=(\mathbb{C}\otimes T\mathcal M')\cap f^{-1}\mathcal N\),
 then one says that
 \(\mathcal{M'}\) is a {\it CR submanifold\/} of \(\mathcal{M}\)
 \cite{Hill99}.
\end{defin}
There is a convenient characterization of CR submanifolds
in terms of an atlas of CR charts:
\begin{prop}\label{p:CRsub}
Let \(f:\mathcal{M'}\to\mathcal{M}\) define \(\mathcal{M'}\)
as a submanifold of the CR manifold \((\mathcal{M},\mathcal{N})\).
Let \eqref{e:CRchart} be a CR chart on \(\mathcal{M}\) and \(\omega\)
the corresponding local section of the canonical bundle.
If, for every such chart,
\begin{equation}\label{e:crsub}
  f^*\omega=0
\end{equation}
and one can find \(n-1\) linear combinations \(({\mu'}^1,\dots,{\mu'}^{n-1})\)
of the forms \((\mu^1,\dots,\mu^{n})\) such that \(\omega'=
f^*(\kappa\wedge{\mu'}^1\wedge \dots\wedge{\mu'}^{n-1})\neq 0\),
then
\begin{equation*}
\mathcal N'=\{w\in\mathbb C\otimes T\mathcal{M'}\mid
w\lrc\omega'=0\}
\end{equation*}
 defines on \(\mathcal{M}'\) the structure of a
CR submanifold of
%the CR manifold
\((\mathcal{M},\mathcal{N})\).
\end{prop}

\begin{proof}
For every \(p\in\mathcal M'\) the monomorphism \(T_p f\), after extension
 to \(\mathbb{C}\otimes T_p\mathcal{M}'\to\mathbb{C}\otimes
  T_{f(p)}\mathcal{M}\),
  restricts
to a monomorphism  \(\mathcal N_p'\to\mathcal N_{f(p)}\) and
the epimorphism \((T_p f)^*\) restricts to an epimorphism
\(\mathcal N^0_{f(p)}\to\mathcal {N'}^0_p\).
Note that \((T_p f)^*(\mathcal N^0_{f(p)}\cap\overline{\mathcal N^0}_{\!\!f(p)}\ts)\)
is the complexification of a real line bundle: it coincides with
\(\mathcal {N'}^0_p\cap\overline{\mathcal {N'}^0_p}\). Therefore, given
a local basis as in \eqref{e:CRchart}, one has
\begin{equation*}(T_p f)^*(\kappa\wedge\mu^1
\wedge\dots\wedge\mu^n)=0
\end{equation*}
and one can choose \(n\) linear combinations of the forms \eqref{e:CRchart}
at \(f(p)\), \(\kappa\) being one of them,
which are mapped by \((T_p f)^*\) to a basis of \(\mathcal {N'}^0_p\).
\end{proof}

\section{hermite and  robinson structures}
\subsection{Almost Hermite and almost Robinson structures}
\begin{defin}
An {\em N-structure} on a Riemannian manifold \((M,g)\)
 of even dimension \(\geqslant 4\), is a complex vector
 subbundle \({N}\) of the
complexified tangent bundle \(\mathbb{C}\otimes TM\) such that,
for every \(p\in M\), the fiber \({N}_p\) is \(\it mtn\).
\end{defin}
It is known that, if \((M,g)\) is proper Riemannian, then
an \(N\)-structure on \(M\) is equivalent to
that of an {\it almost Hermite\/} manifold;
the orthogonal almost complex structure \(J\) on \(M\) is defined as in
\eqref{e:defJ} (see, e.g.,
Ch. IX \S 4 in \cite{KobNom}).

\begin{defin}
An {\it almost Robinson\/} manifold is a Lorentzian
manifold with an \(N\)-structure.
\end{defin}
 In this case,
the intersection \({N}\cap\bar{N}\)
is the complexification of a line bundle \({K}\subset TM\);
its fibers are null; they are tangent to a foliation of
\(M\) by null curves.
An  almost Robinson structure on \(M\) is said to be {\it regular\/}
if the set \(\mathcal M\) of the leaves of the foliation
defined by \(K\) has the structure of a manifold such that
the natural map \(\pi:M\to\mathcal M\) is a submersion.
>From now on, only such regular structures will be considered.

\subsection{The integrability condition}
\begin{defin}
The \(N\)-structure \(N\to M\) on a Riemannian manifold \((M,g)\)
is said to be {\it integrable\/} if
\begin{equation}\label{e:Sec}
  [\Sec N,\Sec N]\subset\Sec N.
\end{equation}
\end{defin}
Dually, the integrability condition is
\begin{equation}\label{e:Sec0}
  \rmd \Sec N^0\subset \Sec N^0\wedge \Sec (\mathbb{C}\otimes
  T^*M).
\end{equation}

In the proper Riemannian case, condition \eqref{e:Sec} is equivalent
to the vanishing of the Nijenhuis (torsion) tensor of the almost complex
structure \(J\) and, by the celebrated Newlander--Nirenberg
theorem, it implies that \(M\) is a Hermite manifold; see Ch. IX \S
2 and 4 in \cite{KobNom}.

%From now on we restrict ourselves to  the Lorentzian case.
\begin{defin}
A {\it Robinson manifold\/} is an  almost Robinson manifold
with an  integrable \(N\)-structure.
\end{defin}

Let \(\omega\) be defined as in \eqref{e:defomega}. It characterizes
\(N\),
\begin{equation}\label{e:Neta}
  N_p=\{w\in \mathbb{C}\otimes T_p M\mid w\lrc \omega=0\}.
\end{equation}
In view of \eqref{e:NN0}, the integrability
condition \eqref{e:Sec0} of Robinson manifolds is
of the same {\it form\/} \eqref{e:int} as for CR structures.

\begin{thm}\label{t:Rman}
Consider a Robinson manifold \(M\) of dimension \(2n+2\).
Let \((\phi_t)\) be the flow generated by a
 vector field \(k:
M\to K\), where \(K\subset TM\) is the null line bundle defined
by \(N\cap \bar{N}=\mathbb{C}\otimes K\), then

\paritem[(i)] the \(N\)-structure on \(M\) is invariant with
respect to the action of the flow \((\phi_t)\) and the trajectories of
\((\phi_t)\)  are null geodesics;

\paritem[(ii)] the \(N\)-structure on \(M\) defines a
 Cauchy--Riemann structure on the quotient manifold \(\mathcal M\);

\paritem[(iii)] the \(2n\)-dimensional
fibers of the bundle \(K^\perp/K\to M\)
have a complex structure and a positive-definite quadratic form, induced
by \(g\).
\end{thm}
\begin{proof}(i) Let \((\kappa,\mu^1,\dots,\mu^n)\) be as
in \eqref{e:NN0}; in view of the reality of \(\kappa\),
the integrability condition \eqref{e:Sec0}
is equivalent to
\begin{equation}\label{e:dlambda}
\rd \kappa=\kappa\wedge\rho+\rmi
\sigma_{\alpha\beta}{\mu}^\alpha\wedge\bar{\mu}^\beta,
\end{equation}
and
\begin{equation}\label{e:dmu}
  \rd \mu^\alpha=\kappa\wedge\varsigma^\alpha+\mu^\beta\wedge{\tau_\beta}^\alpha,
  \;\;\;\alpha=1,\dots,n,
\end{equation}
where
\(\rho, \varsigma^\alpha\) and \({\tau_\beta}^\alpha\)
 are one-forms and the \(\sigma\)s are
functions such that  \(\sigma_{\alpha\beta}=
\overline{\sigma_{\beta\alpha}}\).
 It follows from  \eqref{e:Neta} that the invariance
of \(N\) with respect to \((\phi_t)\) is equivalent to
\(L(k)\omega\Vert\omega\);  this relation  follows from
\eqref{e:dlambda} and \eqref{e:dmu}. Moreover,
 equation \eqref{e:dlambda} implies
\begin{equation}\label{e:lLl}
  \kappa\wedge L(k)\kappa=0.
\end{equation}
In view of  \eqref{e:kkappa} one  can take \(\kappa=g(k)\) so that
  \(L(k)\kappa=(L(k) g)(k)=\nabla_k\kappa\); this shows
 that \eqref{e:lLl} is equivalent to the geodetic condition
\(\nabla_k k\Vert k\).

\noindent (ii) It follows from (i) that the distribution \(N\subset
\mathbb{C}\otimes TM\) descends to a distribution
\(\mathcal{N}\subset\mathbb{C}\otimes T\mathcal M\); its fibers are
of complex dimension \(n\) and
 \(\mathcal{N}\cap\overline{\mathcal{N}}=\pmb 0\).
 Moreover, the integrability of \(N\) implies that
of \(\mathcal{N}\).

\noindent (iii) Only the complex structure requires a construction:
since
 \[K^\perp =\re (N+\bar{N}),\]
  one can put
\(J(w+\bar{w} \bmod K)=\rmi(w-\bar{w}) \bmod K\) for \(w\in N\).
\end{proof}
Note that if \(k\) and \(k'\) are two sections of \(K\to M\),
nowhere vanishing on open subsets \(U\) and \(U'\) of \(M\),
respectively, then \(k'=fk\), where \(f\) is a nowhere zero function
on \(U\cap U'\). If \((\phi_t)\) and \((\phi'_t)\) are the flows
generated by \(k\) and \(k'\), respectively, then, on \(U\cap U'\),
the invariance of \(N\) with respect to \((\phi_t)\) is equivalent
to that with respect to \((\phi'_t)\) and the trajectories of these
two flows coincide.

There is a  local converse to Theorem \ref{t:Rman}. Let
\(\mathcal M\) be a \((2n+1)\)-dimensional CR manifold
characterized  by  differential forms as described in Section~\ref{s:CR}.
Put
\begin{equation}\label{e:defM}
  \pi=\mathrm {pr_1}:M=\mathcal M\times\mathbb{R}\to\mathcal M.
\end{equation}
and denote by \(\kappa\), \(\mu^1\),\dots,\(\mu^n\) the pull-backs by
\(\pi\) to \(M\) of the corresponding forms on \(\mathcal{M}\).
Let \(v\) be the canonical coordinate on \(\mathbb{R}\) and
\(k=\partial_v\in\Sec TM\). The collection of forms
\begin{equation}\label{e:forms}
(\kappa,\rd v,\mu^1,\dots,\mu^1,\bar{\mu}^1,\dots,\bar{\mu}^n)
\end{equation}
is a (local) basis of \(\Sec(\mathbb{C}\otimes T^*M)\);
let
\begin{equation*}
(l,k,\bar{ Z}_1,\dots,\bar{ Z}_n, Z_1,\dots, Z_n)
\end{equation*}
 be the dual basis. We shall construct a Robinson
manifold \((M,g,N)\) so that \eqref{e:NN0} holds. With respect
to the basis \eqref{e:forms}, the metric is
\begin{equation*}
  g=g(l)\kappa+g(k)\rd v+g(\bar{ Z}_\alpha)\mu^\alpha
  +g( Z_\alpha)\bar{\mu}^\alpha.
\end{equation*}
Note that since \(k\in \Sec (N+\bar{N})^\perp\),
one has \(g(k)=g(k,l)\kappa\); therefore \(g(k,l)\neq 0\).
Defining
\(
\lambda=g(l)+ g(k,l)\rd v+g(\bar{ Z}_\alpha,l)\mu^\alpha+
g( Z_\alpha,l)\bar{\mu}^\alpha
\)
so that \(k\lrc \lambda=2g(k,l)\), putting
\(
  g_{\alpha\beta}=2g(\overline{ Z}_\alpha, Z_\beta)=
  \bar{g}_{\beta\alpha}
\),
one obtains
\begin{equation}\label{e:gforRob}
g=\kappa\lambda+g_{\alpha\beta}\mu^\alpha\bar{\mu}^\beta.
\end{equation}
This concludes the proof of
\begin{prop}\label{p:iRman}
Locally, every Robinson \((2n+2)\)-manifold \((M,g,N)\), having
\(\mathcal{M}\) as  the associated
CR manifold, is of the form \eqref{e:defM}
with a metric given by \eqref{e:gforRob}, where the forms
\(\kappa\), \(\mu^1\),\dots,\(\mu^n\) are obtained by pull-back
of the corresponding forms on \(\mathcal{M}\), the functions
\(g_{\alpha\beta}:M\to\mathbb{C}\) are such that, for every
\(p\in M\), the form \(g_{\alpha\beta}(p)z^\alpha
\bar{z}^\beta\) is Hermitean positive-definite,  \(\lambda\) is
any real 1-form on \(M\) such that \(k\lrc\lambda\) is nowhere \(0\)
and \(N^0=\sspan\{\kappa,\mu^1,\dots,\mu^n\}\).

\end{prop}
\subsection{Four-dimensional Robinson manifolds: space-times with
a non-distorting foliation by  null geodesics}\label{ss:Rob}
The case of dimension 4 is well known, but, since it is also
the most important one, it is worth-while to review it briefly
 here. In a sense made precise below, in this case, unlike
 as in higher dimensions, all information about the Robinson
 structure is encoded
 in the properties of the bundle \(K\).

Let \((M,g)\) be a space- and time-oriented Robinson manifold
of dimension 4 with the bundle \(N\to M\) of {\it mtn\/} spaces.
The fibers of the bundle  \(K^\perp/K\to M\)
are two-dimensional  `screen spaces'. According to part (iii)
of  Theorem \ref{t:Rman},
each screen space has a complex structure, which, {\it in this case\/},
is
 equivalent to a conformal
structure and an orientation; this being preserved by the flow is
equivalent to \cite{RT3}
\begin{equation}\label{e:sf}
L(k) g=\rho g +\kappa\otimes\xi +\xi\otimes\kappa
\end{equation}
for some function \(\rho\) and 1-form \(\xi\). Physicists say
that \(k\) generates a shear-free congruence of null geodesics.
The expression
`shear-free' reflects the non-distorting property property of the flow:
it preserves the conformal structure of the screen spaces.
Conversely, given a bundle
\(K\) of null directions, the space and time orientations of \(M\)
induce an orientation in the screen spaces; together with
the induced  Euclidean metric this determines a complex structure \(J\)
in each screen space. This complex structure defines the bundle
\(
N=\{w\in\mathbb C\otimes K^\perp\mid J(w \bmod \mathbb{C}\otimes
K)=\rmi w\bmod\mathbb{C}\otimes K\}
\)
with {\it mtn\/} fibers.
Equation \eqref{e:sf} implies \([\Sec K,\Sec N]\subset \Sec N\);
in dimension 4 this is enough to establish the validity
of \eqref{e:int}.
In view of this, we shall often denote
by \((M,g,K)\) a Robinson space-time determined by the bundle \(K\)
of null lines satisfying \eqref{e:sf}.
 As a consequence of Proposition \ref{p:iRman} one has
\begin{cor}
Let \(\mathcal M\) be a CR  space.
 Put \(M=\mathcal M\times\mathbb{R}\), denote
by \(v\) a coordinate on \(\mathbb{R}\), put \(k=\partial_v\),
\(K=\sspan k\),
pull-back to \(M\) the forms characterizing the CR structure on
\(\mathcal{M}\) to obtain the pair \((\kappa,\mu)\). Let
\(p:M\to \mathbb{R}^+\) and let \(\lambda\) be a 1-form on
\(M\) such that \(k\lrc\lambda\neq 0\). If
\begin{equation}\label{e:g4}
  g=\kappa\lambda+p\mu\bar{\mu},
\end{equation}
then \((M,g,K)\) is a Robinson space-time and every Robinson space-time
can be locally so described, as a {\em lift\/} of \(\mathcal M\).
\end{cor}
\begin{problem}
Characterize the CR spaces that admit lifts to Einstein--Robinson
space-times.
\end{problem}
\begin{thm}
Let \((M,g,K)\) be a Robinson space-time so that
\(g\) is of the form \eqref{e:g4} and the \(N\)-structure
is characterized by \(N^0=\sspan\{\kappa,\mu\}\).
Given a function
\(\rho:M\to\mathbb{R}^+\) and a 1-form \(\xi\) on \(M\) such that
\begin{equation}\label{e:xi}
  k\lrc (\lambda+\xi)\neq 0,
\end{equation}
 define
\begin{equation*}
  g'=\rho (g+\kappa\xi).
\end{equation*}
Then
\paritem[(i)] \((M,g',K)\) is a Robinson manifold,
\paritem[(ii)] if \(F\) satisfies \eqref{e:starF}-\eqref{e:Fnull2}
on \((M,g,K)\), then it also satisfies these equations on \((M,g',K)\).
\label{t:Bat}
\end{thm}
\begin{proof}
(i) One has \(g'=\rho(\kappa\lambda'+p\mu\bar{\mu})\), where \(\lambda'=
\lambda+\xi\) and \(\kappa\wedge\lambda'\wedge\mu\wedge\bar{\mu}\neq 0\)
by virtue of \eqref{e:xi}. Moreover, the bundle \(N\to M\)
does not change under the replacement of \(g\) by \(g'\).

\noindent (ii) The properties
\eqref{e:starF}-\eqref{e:Fnull2} of the form \(F=A\kappa\wedge\mu\)
 also do not change.
\end{proof}
The theorem originates with work of Bateman \cite{Bat10}; see
also \cite{RT3}.
The geometry of \((M,g')\) may be rather different from that of
\((M,g)\); the electromagnetic fields defined by \(F\) in these two
space-times may also be physically distinct. This is illustrated by
the following
\begin{ex} \label{ex:B}
Let \(\mathbb{R}^4\) be the Minkowski space-time.
It is convenient to use
a global coordinate system \((u,v,w)\), where the coordinates
\(u, v\) are real and \(w\) is complex so that
\begin{equation}\label{e:mink}
  g=\rd u\rd v+\rd w\rd\bar{w}.
\end{equation}
Consider the \(N\)-structure
corresponding to \(\sspan\{\rd u,\rd w\}\).
 If \(A(u,w)\) is a function complex-analytic
in \(w\), smoothly depending on \(u\), then the complex 2-form
\begin{equation}\label{e:F=}
  F=A(u,w)\rd u\wedge\rd w
\end{equation}
satisfies equations \eqref{e:starF}-\eqref{e:Fnull2} with \(\kappa=\rd u\);
it describes a {\it plane-fronted\/} electromagnetic wave. If \(A\) depends
on \(u\) only, then \(F\) is a plane wave.

Consider now the open submanifold \(M\) of \(\mathbb{R}^4\) defined by
\(v>0\) and put, for \(m\in\mathbb{R}^+\),
\begin{equation*}
  \rho=v^2(1+\sfrac{1}{4}w\bar{w})^{-2},\quad\rd v+\xi=
  \rho^{-1}(1-2mv^{-1})\rd u+2\rho^{-1}\rd v.
\end{equation*}
Then
\begin{equation*}
  g'=(1-2mv^{-1})\rd u^2+2\rd u\rd v+\rho\rd w\rd\bar{w}
\end{equation*}
and \((M,g')\) describes the Schwarzschild space-time. The form
\eqref{e:F=} corresponds now to a wave with spherical fronts; its
amplitude decreases as \(1/v\) along the null lines of the
expanding congruence generated by \(k=\partial_v\).
\end{ex}

If the CR structure underlying a Robinson space-time \((M,g,K)\) is trivial,
then one can choose coordinates so that \(\kappa=\rd u\)
and \(\mu=\rd w\), as in the last Example.
 In such a case physicists say that \(K\)
defines an \sng\  congruence  {\it without twist\/}.
There are many Einstein--Robinson space-times of this kind. For example,
if the function \(f(u,x,y)\) satisfies the Laplace equation,
\(\partial^{\ts 2}_x f+\partial^{\ts 2}_y f=0\), then
the {\it plane-fronted gravitational wave},
\begin{equation*}
  g=f(u,x,y)\rd u^2 +2\rd u\rd v+\rd x^2 +\rd y^2,
\end{equation*}
has vanishing Ricci tensor, but is not flat unless \(f\) is linear
in \(x\) and \(y\). Its Weyl tensor is of type N.
The plane-fronted waves are among Lorentzian
analogs of K\"ahler manifolds of proper Riemannian geometry: their
bundle \(N\to M\) is invariant with respect to parallel transport.
\begin{problem}
In dimension \(\geqslant 4\), develop a theory of Robinson
manifolds analogous to K\"ahler manifolds.
\end{problem}

`Twisting' congruences, characterized by \(\rd\kappa\wedge\kappa\neq 0\),
are more interesting; the Kerr space-time, describing a black hole
arising from the collapse of a rotating star, is a Robinson manifold
with a twisting congruence.

\begin{ex}\label{ex:Rob}
 In Minkowski space-time, one of the first twisting
shear-free congruences of null lines was described by Robinson around
1963; it played a major role in the emergence of Penrose's
twistors \cite{Pen67,Pen87}.
Robinson established that the metric tensor
\begin{equation}\label{e:Robcongr}
    g=(\rd u+\rmi (z\rd \bar{z}-\bar{z}\rd z))\rd v+(v^2+1)\rd z
\rd\bar{z},\quad z=x+\rmi y
\end{equation}
is flat and the \sng\ congruence generated by \(\partial_v\)
is twisting.
The complex 2-form
\(F=A(x,y,u,v)\kappa\wedge
(\rd x+\rmi \rd y)\)
 is self-dual and Maxwell's equations \(\rd F=0\)
reduce to \(\partial A/\partial v=0\) and the equation
\({ Z}\lrc\rd A=0\),
where
\(
 { Z}=\partial_x+\rmi\partial_y
  -\rmi(x+\rmi y)\partial_u
\)
is an  operator on \(\mathbb{R}^3\) introduced by Hans Lewy in 1957. He
constructed a smooth function \(h\) such that the equation
\({ Z}\lrc\rd A=h\) has no solution, even locally.

 The underlying CR geometry on \(\mathcal M=\mathbb{R}^3\) with
coordinates \(u,z=x+\rmi y\) is given by the pair
\((\kappa=\rd u+\rmi (z\rd \bar{z}-\bar{z}\rd z),\mu=\rd x+\rmi\rd y)\).
Two solutions of \eqref{e:embed}
are \(z_1=x+\rmi y\) and \(z_2=u+\half \rmi (x^2+y^2)\) so that
equation \eqref{e:G=0} is now that of the hyperquadric,
\(\rmi(\bar{z}_2-z_2)-\vert z_1\vert^2=0\).
The biholomorphic map
\begin{equation*}
w_1=\sqrt 2\frac{z_1}{z_2+\rmi},\quad w_2=\frac{z_2-\rmi}{z_2+\rmi}
\end{equation*}
transforms the hyperquadric into the 3-sphere of equation
 \[\vert
w_1\vert^2+\vert w_2\vert^2=1.\]
This is the most symmetric,
non-trivial, 3-dimensional CR geometry: its group of automorphisms
is \(\SU_{2,1}\). The CR structure on \(\mathbb{S}_3\) can be
viewed as obtained from the complex structure of
\(\mathbb{S}_2=\CP_1\) via the Hopf map.

 Several solutions of Einstein's equations
 admit this congruence. As an example, we show this for the
 {\it G\"odel universe} \cite{Koch}. Take its metric in the form
 given in \cite{RT3},
\begin{equation*}
(\rd X^2+\rd Y^2 -2(Y\rd U-\rd X)(Y\rd V-\rd X))/Y^2.
\end{equation*}
Its Weyl tensor is of type D: the null vector fields \(k=\partial_V\)
and \(l=\partial_U\) generate each an \sng\  congruence.
Consider \(k\); the corresponding CR structure on \(\mathbb{R}^3\) with
coordinates \((U,X,Y)\) is given by \(\kappa=\rd X-Y\rd U\)
and \(\mu=\rd X+\rmi\rd Y\). Introduce new local coordinates \((u,x,y)\)
in \(\mathbb{R}^3\) by
\(
  u=X,\quad z=x+\rmi y=\sqrt{Y}\exp(-\half\rmi U)
\).
One then obtains
\(
  \kappa=\kappa',\quad \mu=\kappa'+2\rmi\bar{z}\mu',
\)
where
\begin{equation*}
  \kappa'=\rd u+\rmi(z\rd\bar{z}-\bar{z}\rd z),\quad \mu'=\rd z.
\end{equation*}
The pair \((\kappa',\mu')\) defines the same CR structure as the pair
\((\kappa,\mu)\): it is  that of the hyperquadric.

\end{ex}
\subsection{The Goldberg--Sachs theorem}\label{s:GS}
Consider a 4-manifold \((M,g)\) that is either proper Riemannian
or Lorentzian. An \(N\)-structure on \(M\) can be (locally) given
by a field \(\varphi\) of chiral spinors: one uses `point by
point' the definition \eqref{e:Nvar}.
\begin{thm}\label{t:GS}
{\rm (i)} If the \(N\)-structure \(N(\varphi)\) is integrable,
then the chiral spinor \(\varphi\) is an eigenspinor of the Weyl
tensor.

\noindent{\rm (ii)} If \((M,g)\) is conformal to an Einstein manifold,
then \(N(\varphi)\) is integrable if, and only if, the chiral
spinor field \(\varphi\) is a repeated eigenspinor of the Weyl
tensor.
\end{thm}
For space-times, the theorem was established by Goldberg
and Sachs \cite{GoldSachs62}. Its extension to the proper Riemannian
case is due to Pleba\'nski, Hacyan, Przanowski and Broda \cite{PlebHac75,PB}.

\begin{problem}
Find a generalization of the Goldberg--Sachs theorem to
manifolds of dimension \(>4\).
\end{problem}

In the Lorentzian case, it follows from Theorem \ref{t:GS}
and the algebraic classification
of Weyl tensors that a space-time which is  conformally Einstein, but
  not conformally flat, can have at most 2 distinct \sng\ congruences
   (type D). The following example shows that there are
  non-conformally flat space-times admitting 3 such distinct congruences;
  we do not know whether there are  space-times with
   \(\mathsf C\neq 0\) and 4 distinct
  congruences of this type.

\begin{ex}
Consider a space-time \(M=\mathbb{R}^4\) with the real coordinates
\(u,v\) and a complex coordinate \(w\). Let the metric tensor be
\(g=\lambda\kappa+\mu\bar{\mu}\),
where
\begin{equation*}
 \kappa=\rd u+\half\rmi(w\rd\bar{w}-\bar{w}\rd w),\quad
\lambda=\rd v-\half\rmi(w\rd\bar{w}-\bar{w}\rd w),\quad
\mu=(w+\bar{w})\rd w.
\end{equation*}
This space-time admits three congruences of shear-free
null geodesics : those generated by the vector fields \(k_1=\partial_u\) and
\(k_2=\partial_v\) are twisting and are both
equivalent to the Robinson congruence. The congruence generated
by
\[k_3=\partial_v-\partial_u+
2\rmi (w+\bar{w})^{-1}(\partial_{\bar{w}}-\partial_{w})\]
is \sng\ and  has vanishing twist. The space-time
\((M,g)\) has a Weyl tensor of type I and does not admit any
other \sng\ congruences.
\end{ex}

 \subsection{Remarks on the embeddability problem}
The property  of a CR space \(\mathcal{M}\) to be embeddable is relevant to
the local existence of a non-zero, null solution of Maxwell's equations on
 space-times obtained as lifts of \(\mathcal{M}\).
If \(\mathcal{M}\) is embeddable, if the forms \(\kappa\) and \(\mu\)
are as in \eqref{e:emb}, and \(g\) is given by \eqref{e:g4}, then
\(F=A(z_1,z_2)\kappa\wedge\mu\) satisfies \eqref{e:starF}-\eqref{e:Fnull2}
for every function \(A\) holomorphic in its two arguments.
In fact, less is required for the local existence of such an \(F\):
if the canonical bundle of \(\mathcal{M}\) admits a locally defined
closed section \(\omega\), then its pull-back to \(M\) can
be taken as \(F\).

It is now known that there are CR  spaces
that are non-embeddable, but have {\it one\/}
 solution of \eqref{e:embed} \cite{Rosay89};
by the results of \cite{Taf85}, extended to higher dimensions in
\cite{Jac87}, such CR spaces do not admit closed, non-zero sections of
their canonical bundle.
Therefore, space-times constructed  as lifts of these CR
spaces  do not admit any associated non-zero null solutions  of
Maxwell's equations.
There are examples of non-embeddable 7-dimensional CR manifolds that
have non-zero, closed, sections of their canonical bundle, but it is
not clear whether there are such examples in dimensions 3 and 5. Further
remarks on this subject are in \cite{T99forES}.

Lewandowski, Nurowski and Tafel \cite{LNT90} established the following
\begin{thm}
If the CR space \(\mathcal{M}\) lifts to an Einstein--Robinson
space-time, then \(\mathcal{M}\) is locally embeddable.
\end{thm}

\section{the kerr theorem}\label{sss:Kerr}
The Kerr theorem provides a method for constructing all integrable
analytic \(N\)-structures in Minkowski space-time  \((M,g)\);
even though it is well-known, we present it here because of its importance.
 See  \cite{PenMcC73,PenRin,Taf85} for further details and references.
Consider the coordinate system and metric \eqref{e:mink}
as given in Example \ref{ex:B}. The manifold of
 all {\it mtn\/} subspaces of one chirality of the complexified Minkowski space
 \(\mathbb{C}^4\)
 is \(\SO_4/\U_2=\CP_1\).

Let \(z\in\mathbb C\) and define
\begin{subequations}\label{e:forKerr}
\begin{align}
k_z=\partial_v-z\partial_w-\bar{z}\partial_{\bar{w}}-
  z\bar{z}\partial_u,\\
  \kappa_z= \rd u-z\rd\bar{w}-\bar{z}\rd w-
 z\bar{z}\rd v,\\
\mu_z=\rd w +z\rd v,\quad\text{and}\quad\lambda_z=\rd v.
\end{align}
\end{subequations}
The map \((\kappa_0,\mu_0,\lambda_0)\mapsto (\kappa_z,\mu_z,\lambda_z)\) is
a proper Lorentz transformation. It is induced by the homomorphisms
\(\mathbb{C}\to\SL_2(\mathbb{C})\to\SO_{3,1}\). The pair
\((\kappa_z,\mu_z)\) defines an {\it mtn\/} subspace \(N_z\)
such that \(\re (N_z\cap\bar{N}_z)=\dir k_z\). The subspace corresponding
to the `point at infinity' of \(\CP_1=\mathbb{C}\cup\{\infty\}\) is
defined by the pair \((\rd v,\rd\bar{w})\) and \(k_\infty=\partial_u\).
Assume now \(z\) to be a complex {\it function\/} on \(M\)
such that its real and imaginary parts are real-analytic
functions of the coordinates \(u,v\), \(\re w\) and \(\im w\).
At every point \(p\) of \(M\)
the pair  \((\kappa_{z(p)},\mu_{z(p)})\) defines an {\it mtn\/}
subspace  of \(\mathbb{C}\otimes T_p M\).
According to \eqref{e:int}, the \(N\)-structure defined by
 \((\kappa_z,\mu_z)\)  is integrable if, and only if,
\begin{equation}\label{e:lammu}
  \rd \kappa_z\wedge\kappa_z\wedge\mu_z=0\;\;\;\text{and}\;\;\;
\rd \mu_z\wedge\kappa_z\wedge\mu_z=0.
\end{equation}
A simple calculation shows that equations \eqref{e:lammu}
reduce to
\begin{align*}
\rd v\wedge\rd z\wedge\rmd  (u-z\bar{w})\wedge\rmd (w+zv)=0,\;\\
\rd \bar{w}\wedge\rd z\wedge\rmd (u-z\bar{w})\wedge\rmd  (w+zv)=0,
\end{align*}
and are thus equivalent to
\begin{equation}\label{e:dmulam}
\rmd (u-z\bar{w})\wedge\rmd  (w+zv)\wedge\rd z=0.
\end{equation}

By the implicit function theorem,  equation \eqref{e:dmulam} implies,
locally, the existence of  a holomorphic
function \(H(z_1,z_2,z_3)\) of three complex variables such that
\begin{equation}\label{e:Kerr}
  H(u-z\bar{w},w+zv,z)=0.
\end{equation}

This proves a theorem attributed to Kerr:
\begin{thm}
 Locally, every integrable
analytic \(N\)-structure
in Min\-kow\-ski space-time \(\mathbb{R}^4\)
is given either by the pair \((\rd v,\rd\bar{w}) \) or by \eqref{e:forKerr},
where \(z:\mathbb{R}^4\to\mathbb{C}\) is a solution of \eqref{e:Kerr}
and \(H\) is a holomorphic
function of three complex variables such that \(\rd H\neq 0\).
\end{thm}
Denoting \(H_1=\partial H/\partial z_1\), etc.,
 one obtains by differentiation of \eqref{e:Kerr}
\begin{equation*}%\label{e:dH}
  H_1\kappa_z+(H_2+\bar{z}H_1)\mu_z+(H_3-\bar{w}H_1+vH_2)\rd z=0.
\end{equation*}
The condition \(\rd H\neq 0\) implies \(H_3-\bar{w}H_1+vH_2\neq 0\).
If \(H_1=H_2=0\), then \(z=\mathrm{const.}
\) and the \(N\)-structure is trivial, i.e. reducible, by a Lorentz
transformation of the coordinates,  to
\(\kappa_0 =\rd u\) and \(\mu_0=\rd w\).
Define
\begin{equation}\label{e:defuz}
    u_z=u-z\bar{w}-\bar{z}w-z\bar{z}v\quad\text{and}\quad
w_z=w+zv.
\end{equation}
Since
\begin{equation}\label{e:Luw}
    L({k_z})u_z=0\quad\text{and}\quad L({k_z})w_z=0,
\end{equation}
 the functions \(u_z\) and \(w_z\)
descend to the CR manifold \(\mathcal M\)
obtained from \(M\) as described in Theorem \ref{t:Rman}.
Moreover, the pair \((\kappa_z,\rd w_z)\)
defines the same \(N\)-structure  on \(M\)
 as the pair \((\kappa_z,\mu_z)\).
The pair \((\kappa_z,\rd w_z)\) defines the CR structure on
\(\mathcal M\).

Assume now that \(H_1\) and/or \(H_2\neq 0\). Equation \eqref{e:Kerr}
can be written as
\begin{equation*}%\label{e:Hnew}
    H(u_z+\bar{z}w_z,w_z,z)=0
\end{equation*}
and shows that \(w_z\)  is a function of \(z,\bar{z}\) and \(u_z\) only.
The integrability condition \(\rd \kappa_z\wedge\kappa_z\wedge
\mu_z=0\) is now satisfied identically and \(\rd \mu_z\wedge
\kappa_z\wedge\mu_z=0\) is equivalent to
\begin{equation}\label{e:CReq}
    \frac{\partial w_z}{\partial{\bar{z}}}-w_z\frac{\partial w_z}
    {\partial{u_z}}=0.
\end{equation}
Using \eqref{e:CReq} one obtains
\begin{equation*}
    \rd w_z=\frac{\partial w_z}{\partial{u_z}}\kappa_z+
    (\frac{\partial w_z}{\partial z} -
\bar{w}_z\frac{\partial w_z}{\partial{u_z}})\rd z.
\end{equation*}
This shows that the pair \((\kappa_z,\rd z)\) defines on
\(\mathcal M\) the same CR structure as the pair \((\kappa_z,\rd w_z)\).
Let \((\partial_{u_z},\bar{ Z}, Z)\) be
the frame on \(\mathcal M\) dual
to the coframe \((\kappa_z,\rd z,\rd\bar{z})\) so that
\begin{equation*}
 Z=\frac{\partial}{\partial{\bar{z}}}-w_z\frac{\partial}{\partial{u_z}}.
\end{equation*}
Equation
\eqref{e:CReq} is now interpreted as a tangential Cauchy--Riemann  equation,
 \( Z\lrc\rd w_z=0\).

The map \((u,v,w)\mapsto (u_z,v,z)\) is a local
diffeomorphism.
 This is seen by computing the
volume form on \(M\),
\begin{equation*}
    \rmi \rd u\wedge\rd v\wedge\rd w\wedge\rd\bar{w}=
\rmi \vert \bar{ Z}\lrc\rd w_z-v\vert^2\rd u_z\wedge\rd v\wedge
\rd z\wedge\rd\bar{z},
\end{equation*}
where use has been made of \eqref{e:CReq}. The distribution
\(\ker \kappa_z\) is integrable
if, and only if,
 \(\bar{ Z}\lrc\rd w_z\) is real. Dropping the subscripts
 \(z\), one has
\begin{cor}
Let \((u,v,z)\) be a local coordinate system on \(M\), let
\(w(u,z,\bar{z})\) be a smooth, complex-valued function satisfying
\begin{equation*}
    \partial_{\bar{z}}w-w\partial_u w=0
\end{equation*}
and put \(\kappa=\rd u+\bar{w}\rd z+w\rd\bar{z}\),
\(\mu=\rd w-v\rd z\). The metric
\begin{equation}\label{e:metr}
    g=\kappa\rd v+\mu\bar{\mu}
\end{equation}
is flat and the vector field \(k=\partial_v\) generates an
expanding \((\mathrm{div}\ts k\neq 0)\) \sng\ congruence.
\end{cor}
\begin{ex}
If \(w=\rmi z\), then \eqref{e:metr}  assumes
the form \eqref{e:Robcongr}
and corresponds to the Robinson congruence of Example \ref{ex:Rob}.

\end{ex}

\section{twistor bundles}%\label{s:tw}
Recall a general idea in geometry: if one  wishes to study a structure,
but there is no distinguished structure, then it is appropriate to consider
 the set of all such structures.

Given an oriented Riemannian \(2n\)-manifold \((M,g)\)
(conformal
geometry suffices), define
its {\it twistor bundles\/} \(P_\pm\) to have, as the total sets,
the collections of all {\it mtn\/} subspaces of \(\mathbb{C}\otimes TM\) of the
\(\pm\) chiralities. These are bundles with fiber \(\SO_{2n}/\U_n\),
which has a canonical metric and complex structure.
If \(\star^2=-\id\), then complex conjugation in \(\mathbb{C}
\otimes TM\) changes the chirality of the {\it mtn\/} subspaces;
this induces an isomorphism of the bundles \(P_+\) and \(P_-\).
They are then identified and denoted by \(P\): such is the case when
\((M,g)\) is a space-time.
The Levi-Civita connection on \(M\) induces a horizontal distribution
on \(P_\pm\); together with the canonical metric on the
fibers, this defines a metric  and
a canonical \(N\)-structure on \(P_\pm\), which need not
be integrable. If \((M,g)\) is proper Riemannian (resp., Lorentzian),
then so is \(P_\pm\) and its canonical \(N\)-structure defines on \(P_\pm\)
the structure of an almost Hermite (resp., almost Robinson) manifold.
\begin{thm} If \(M\) is a space-time, then the
integrability of the canonical \(N\)-structure on its twistor bundle
\(P\) is equivalent to \(\mathsf C=0\).  If
\(M\) is a \(4\)-dimensional proper Riemannian manifold,
 then the canonical \(N\)-structure
on \(P_\pm\) is integrable if, and only if, \(\mathsf C_\pm=0\).
\label{t:AtH78}
\end{thm}
In the Lorentzian case, the theorem was established by Penrose
in the course of work that led to his fundamental
twistor programme; see  \cite{PenRin}
and the references given there.
The proof in the proper Riemannian case is due to Atiyah, Hitchin and Singer
\cite{AtH78}.

\subsection{The Kerr theorem revisited}\label{ss:Kerrrev}
Let \((M=\mathbb{R}^4,g)\) be the Minkowski space-time. According to
Theorem \ref{t:AtH78}, its twistor bundle \(P\) is a Robinson manifold
so that there is the associated 5-dimensional CR manifold
\(\mathcal{P}\).
The twistor bundle \(P\) is identified with the set of
  null directions in the tangent spaces at all points of \(M\).
  Its typical fiber is the `celestial sphere'
   \(\mathbb{S}_2\approx\CP_1\) so that \(P=M\times \CP_1\).
 Locally, the bundle \(P\to M\) can be conveniently described
  as follows. Let \((u,v,w)\) be a coordinate system on
  \(M\), as in \eqref{e:mink}. A number \(z\in\mathbb{C}\) defines
  a null direction \(\dir k_z\) at \((u,v,w)\), parallel to the
  vector \(k_z\) given in \eqref{e:forKerr}.
 A point of \(P\) is given by the sequence
  \((u,v,w,\dir k_z)\) or, equivalently, by the sequence
\((u,v,w,z)\), i.e. by a sequence of 6 real functions; they
 provide a
 convenient coordinate system on
 \(P\). In these coordinates, the metric tensor on \(P\) is
\(
  \rd u\rd v+\rd w\rd \bar{w} +(1+\sfrac{1}{4}z\bar{z})^{-2}\rd z\rd\bar{z}
\).
The canonical  \(N\)-structure
 on \(P\) is given by
\(
  N_P^0=\sspan\{\kappa_z,\mu_z,\rd z\}
\).
Its integrability  is easily checked
by computing \(\omega_z=\kappa_z\wedge\mu_z\wedge\rd z\) and verifying that
equations \eqref{e:int} are satisfied. The line bundle \(N_P\cap\bar{N}_P\to P\)
is spanned by \(\dir k_z\).

Consider now the  CR manifold \(\mathcal{P}\) associated with \(P\)
as in Theorem~\ref{t:Rman} and the functions defined in \eqref{e:defuz}.
In view of \eqref{e:Luw} and
\(L({k_z}) z=0 \), the sequence \((u_z,w_z,
z)\) of 5 real functions
descends to \(\mathcal{P}\) and provides a coordinate system
on that manifold. Its CR structure
is  embeddable: three solutions of \eqref{e:embed} are
\(z_1=u-z\bar{w}\),  \(z_2=w+zv\) and \(z_3=z\).
Consider a regular congruence \(K\) of null
lines on \(M\) {\it which need not be shear-free\/}.
The set \(\mathcal{M}\) of these lines is a
3-dimensional manifold.
 There is the
map \(f:\mathcal{M}\to\mathcal{P}\) that sends an element
of the congruence on \(M\) to its lift to \(\mathcal{P}\),
\begin{equation*}
\begin{CD}
P@>{\mathrm{can}}>>\mathcal P\\
@V{\mathrm{can}}VV @AA{f}A\\
M@>>\pi>\mathcal M
\end{CD}
\end{equation*}
\begin{thm}
The congruence \(K\) of null lines on Minkowski space-time
is shear-free if, and only if, the map \(f:\mathcal{M}\to\mathcal{P}\)
defines on \(\mathcal{M}\) the structure of a CR submanifold of \(\mathcal{P}\).
\end{thm}
\begin{proof}
Let \(z:M\to\mathbb{C}\) be the function defining
the congruence \(K\) of null lines. The map \(f\circ \pi:M\to\mathcal P\)
sends \((u,v,w)\) to \((u_z,w_z,z)\) with \(z\) evaluated at \((u,v,w)\).
A section of the canonical bundle of the CR manifold \(\mathcal P\) is
\(\omega=\rmd(u-z\bar{w})\wedge\rd w_z\wedge\rd z\). According to
\eqref{e:dmulam},  the pull-back
\((f\circ\pi)^*\omega\) vanishes if, and only if, the null
geodetic congruence \(K\) is shear-free. Since \(\pi\) is a surjective
submersion, this holds only whenever \eqref{e:crsub} is satisfied.
\end{proof}

The image of \(\mathcal{P}\) in \(\mathbb{C}^3\) is the
hypersurface (`generalized hyperquadric')  of equation
\begin{equation}\label{e:crmin}
z_3-\bar{z}_3+z_1 \bar{z}_2 -\bar{z}_1z_2=0.
\end{equation}
Every point of this hypersurface corresponds to a null line  \(
l:\mathbb{R}\to M\) given,
in the coordinate system \((u,v,w)\) on \(M\),  by
\begin{equation*}
l(t)= (\half(z_3+\bar{z}_3+z_1 \bar{z}_2 +\bar{z}_1z_2)
-z_1\bar{z}_1 t,\,t,\, z_2-z_1 t)
\end{equation*}
so that \(l(v)=(u,v,w)\) and \(\rd l/\rd t=k_z\).
All null lines in \(M\), except those parallel to \(\partial_u\),
can be  obtained by this `Penrose correspondence'
 between \(M\) and \(\mathcal{P}\).
Consider now the embedding
\begin{equation*}
 f: \mathbb{C}^3\to\CP_3,\quad f(z_1,z_2,z_3)=\dir (1+\rmi z_3,
  z_1-\rmi z_2,1-\rmi z_3,z_1+\rmi z_2).
\end{equation*}
The image of \(\mathbb{C}^3\) by \(f\) is \(\CP_3\) with
a \(\CP_2\) removed.
 The image of the hypersurface \eqref{e:crmin} by \(f\) is an open
and dense submanifold of  the manifold \(\mathcal{P}_0\)
of {\it null twistor\/} directions
\begin{equation}\label{e:nt}
  \{\dir(w_1,w_2,w_3,w_4)\in\CP_3\mid \vert w_1\vert^2+\vert w_2\vert^2
  -\vert w_3\vert^2-\vert w_4\vert^2=0\}.
\end{equation}

Penrose \cite{PenMcC73} proved the following fundamental
\begin{thm}\label{t:pkt}
 If \(M=(\mathbb S_1\times\mathbb S_3)/\mathbb Z_2\)
  is the conformally compactified Minkowski space-time, then
\(P=\CP_3\).
Every analytic CR \(3\)-manifold, defining a Robinson structure in
\(M\), is obtained as the intersection of the \(5\)-dimensional CR
manifold of projective null twistors  \eqref{e:nt} with a complex
analytic \(2\)-dimensional submanifold of \(\CP_3\).
\end{thm}
According to Penrose, a  non-analytic, shear-free and twisting congruence of
  null geodesics in (compactified) Minkowski space-time
can be described as corresponding to a complex surface
\(\varSigma\) in \(\CP_3\)  that
`touches only one side' of the manifold of projective null twistors
\(\mathcal P_0\)
so that the real dimension of  \(\mathcal{P}_0\cap\varSigma\)
 is 3, but the surface cannot be holomorphically extended
  to the other side of \(\mathcal P_0\), see pp.~220--222 in
  \cite{PenRin}.

\subsection{The Kerr theorem in the proper Riemannian setting}
There is an  analog of the Kerr theorem for proper Riemannian
self-dual (or anti-self-dual) 4-manifolds. We only sketch the
idea of the theorem in the {\it local\/}
 setting. According to Theorem \ref{t:AtH78},
the twistor bundle \(P_+\) of such  a self-dual manifold has a canonical
integrable \(N\)-structure defining there the structure of a complex
3-manifold so that there is the fibration
\(\CP_1\to P_+\stackrel{\pi}{\to} M\). Let \(U\) be an
open subset of \(M\) and \(s:U\to P_+\) a local section
of \(\pi\) such that \(s(U)\) is a complex submanifold of \(P_+\).
The restriction of \(\pi\) to \(s(U)\) induces on \(U\) the structure
of a Hermite manifold and all local Hermite structures on \(M\)
can be so obtained. The insistence on locality is essential: for example,
the 4-sphere  has no global complex structure, but it has local
Hermite structures.

 \section*{acknowledgments}
 We thank C.D. Hill, P. Kobak, L.J.  Mason, J. Tafel
  and K.P. Tod for enlightening
 discussions on several topics presented here.
We are grateful to I.~Bia{\l}ynicki-Birula for
the references \cite{Silberstein1907} and \cite{Weber1910}.

Work on this paper was supported in part by the Polish Committee
for Scientific Research (KBN) under grant no. 2 P03B 060 17.

\providecommand{\bysame}{\leavevmode\hbox to3em{\hrulefill}\thinspace}

\end{document}